\documentclass[11pt]{article}
\usepackage{amsmath, amssymb,amsfonts,amsthm}

\newtheorem*{remark}{Remark}
\newtheorem*{theoremA}{Theorem A}
\newtheorem*{theoremB}{Theorem B}
\newtheorem*{gracies}{Acknowledgements}

\newcommand{\RN}{\mathbb{R}^N}

\newcommand{\ep}{\varepsilon}
\newcommand{\ca}{\operatorname{Cap}}

\title{Capacitary differentiability of potentials of finite Radon measures}

\author{Joan Verdera}
\date{}

\begin{document}

\maketitle

\footnotetext{\\*[-11pt]
}

\begin{abstract}
We study differentiability properties of a potential of the type $K\star \mu$, where $\mu$ is a finite Radon measure
in $\RN$ and the kernel $K$ satisfies $|\nabla^j K(x)| \le C\, |x|^{-(N-1+j)}, \quad j=0,1,2.$
 We introduce a notion of differentiability in the capacity sense, where capacity is
classical capacity in the de la Vall\'ee Poussin sense associated with the kernel $|x|^{-(N-1)}.$
We require that the first order remainder at a point is small when measured by means of a normalized weak
capacity ``norm'' in balls of small radii centered at the point. This implies weak $L^{N/(N-1)}$ differentiability and thus
$L^{p}$ differentiability in the Calder\'on--Zygmund sense for $1\le p < N/(N-1)$. We show that $K\star \mu$ is a.e. 
differentiable in the capacity sense, thus strengthening a recent result by Ambrosio, Ponce and Rodiac. We also present an
alternative proof of a quantitative theorem of the authors just mentioned, giving pointwise Lipschitz estimates for $K\star \mu.$
As an application, we study level sets of newtonian potentials of finite Radon measures.
\end{abstract}

\section{Introduction}
Let  $K$ be a kernel twice continuously differentiable off the origin of $\RN$ satisfying  
\begin{equation}\label{ker}
|\nabla^j K(x)| \le C\, \frac{1}{|x|^{(N-1+j)}}, \quad j=0,1,2.
\end{equation}
Given a finite Radon measure $\mu$ we would like to understand the differentiability properties of the locally integrable function
$K\star \mu.$ To obtain sharp results we resort to a notion of differentiability in the capacity sense that was introduced recently
in \cite{CV} to study  differentiability properties of newtonian potentials of finite Radon measures  (logarithmic potentials of 
compactly supported Radon measures in the plane). 
The capacity we use is the classical one, in the de la Vall\'ee Poussin sense,  associated with the kernel $|x|^{-(N-1)}.$ 
 Given a compact $E \subset \RN$ the capacity of $E$ is defined as
\begin{equation*}
 \ca(E) = \sup \{\mu(E) : \mu \ge 0,\;\;\text{spt} \mu \subset E\; \; \text{and}\;\;
 \left(\frac{1}{|x|^{N-1}}\star\mu \right) (x) \le 1, \;\; x \in \RN \}.
\end{equation*}
The capacity of an arbitrary  set is then defined as the supremum of the capacities of compact subsets.
The homogeneity of $\ca$ is $N-1,$ so that the capacity of the ball of center $x$ and radius $r$ is $\ca(B(x,r))= c\, r^{N-1}$ for some constant $c.$

For each finite Radon measure $\mu$ the potential   $u = K \star \mu$  is well defined except for a set of zero capacity. We say that $u$ is  differentiable at the point $a$ in the capacity sense if  there exist a vector $\vec{v} \in \RN,$ called the gradient of $u$, such that
\begin{equation}\label{difcap}
 \lim_{r \rightarrow 0} \frac{1}{\ca(B(a,r))} \sup_{t>0} t \, \ca\left(\{x\in B(a,r) : 
 \frac{|u(x)-u(a) -\vec{v}\cdot(x-a) |}{|x-a|} > t \}\right) =0.
\end{equation}
As we will discuss later  this implies differentiability in the $L^{\frac{N}{N-1}, \infty}$ sense, which in turn yields differentiability
in the $L^p$ sense, $1 \le p < \frac{N}{N-1},$ and thus approximate differentiability. 

\begin{theoremA}\label{TA}
For each kernel $K$ satisfying \eqref{ker} and for each finite Radon measure $\mu$ the potential  $K\star \mu$ is differentiable in the capacity sense
 at almost all points of $\RN.$
\end{theoremA}

Hajlasz studied in \cite{H} the case in which $K \in C^2(\RN\setminus \{0 \})$ is homogeneous of degree $-(N-1)$ and proved
that $K\star \mu$ is approximately differentiable almost everywhere.  Alberti, Bianchini and Crippa \cite{ABC}, working under the 
same hypothesis, showed that $K\star \mu$ is differentiable almost everywhere in the $L^p$ sense, for $1 \le p < \frac{N}{N-1},$ 
which implies approximate differentiability. In a recent preprint Ambrosio, Ponce and Rodiac \cite{APR} have introduced the more 
general class of kernels \eqref{ker} and proved differentiability almost everywhere in the $L^{\frac{N}{N-1}, \infty}$ sense. 
The proofs of the results mentioned above are based on an application of the Calder\'on-Zygmund decomposition of a finite Radon measure. 
We avoid appealing directly to such a decomposition, which is used indirectly via the classical weak type estimates of
Calder\'on -Zygmund operators.

As we mentioned before capacitary differentiability \eqref{difcap} is stronger than $L^{\frac{N}{N-1}, \infty}$ differentiability. 
This follows from the weak type capacitary inequality
\begin{equation}\label{wcap}
 t \, \ca(\{x : \left| \left(\frac{1}{|x|^{N-1}}\star\mu\right)  (x) \right | > t \}) \le \|\mu\|, \quad t >0,
\end{equation}
combined with the elementary estimate
\begin{equation}\label{mecap}
 |E|^{\frac{N-1}{N}} \le C \, \ca(E), \quad E \quad \text{measurable},
\end{equation}
where bars denote Lebesgue measure in $\RN.$

Next result is a pointwise Lipschitz estimate for the potential $K\star \mu,$ which yields immediately approximate differentiability.
A precursor of that is Lemma 9 in \cite{H}, where the Lipschitz estimate is given  for $K\star \mu^b,$  where $\mu^b$ is the bad
part  in the Calder\'on -Zygmund 
decomposition of $\mu$ at a given height $t.$ In \cite{APR} the inequality is proven for the whole potential $K\star \mu$ 
(and for the class of kernels satisfying 
\eqref{ker}), with a fixed dominating function in $L^{1,\infty}(\RN),$ independent of $t.$
We have denoted by $L^{1,\infty}(\RN)$ 
the space of weak $L^1$ functions, that is, the space of  measurable functions $I$ on $\RN$ satisfying
$\sup_{t>0} t \left|\{x \in \RN: |I(x)| >t \}\right| < \infty.$

There are two classical operators that send finite Radon measures into $L^{1,\infty}(\RN),$ with bounds. 
The first is the  Hardy-Littelwood maximal 
operator $M,$ which on a finite Radon measure $\mu$ is defined as 
$$
M(\mu)(x) = \sup_{r>0} \frac{|\mu(B(x,r))|}{|B(x,r)|}, \quad x \in \RN,
$$
where $B(x,r)$ is the ball centered at $x$ of radius $r.$
The second is the maximal singular integral. It turns out that $\nabla K$ is a vector valued Calder\'on -Zygmund kernel 
(see section 2 for more details). The maximal singular 
integral associated with that kernel is
$$
T^*(\mu)(x) = \sup_{\varepsilon >0} \left|\int_{|x-y|>\varepsilon} \nabla K(x-y)\,d\mu(y)\right|, \quad x \in \RN.
$$

\begin{theoremB}\label{TB}
For each kernel $K$ satisfying \eqref{ker} there exists a positive constant $C$ such that for each finite Radon measure $\mu$ 
one can find a non-negative function 
$I \in L^{1, \infty}(\RN)$ satisfying 
 \begin{equation}\label{klip}
 \left| \left(K\star \mu \right)(x)-\left(K\star \mu \right)(y)\right| \le |x-y|\, \left(I(x)+I(y)\right), 
 \quad \text{for almost all}\quad x, y \in \RN,
\end{equation}
and
$$
\sup_{t>0} t \left|\{x : I(x) >t \}\right| \le C \, \| \mu\|.
$$

Indeed, on can take $I= C\,\left( M(\mu) + T^*(\mu)\right)$ in \eqref{klip}.
 \end{theoremB}

 The above result has been recently proven in \cite{APR} except for the last assertion concerning the precise form of $I.$ 
 Our contribution here is to provide a new proof, which does not appeal to a
Calder\'on-Zygmund decomposition of $\mu$ and gives a canonical choice for the dominating function $I.$

In section 2 we prove Theorem A and in section 3 Theorem B. Section 4 is devoted to an application, devised in \cite{APR}, to level sets 
of potentials of the form  $\Phi \star \mu,$ where $\Phi$ is the fundamental solution of the laplacian in $\RN.$

Our notation and terminology are standard. The letter $C$ denotes a constant which may change at each occurence and
that is independent of the relevant parameters under discussion. 

\section{Proof of Theorem A}
The kernel $\nabla K$ has the growth and smoothness conditions of a Calder\'on-Zygmund kernel and it has also the cancellation property
$$
\left| \int_{R_1 < |x| < R_2} \nabla K(x) \,dy \right| = \left|\int_{|x|=R_2} K(x) \frac{x}{|x|}\,d\sigma(x) - 
\int_{|x|=R_1} K(x) \frac{x}{|x|}\,d\sigma(x) \right| \le C.
$$
Therefore the general theory of Calder\'on-Zygmund operators, as presented in Grafakos' book \cite[Chapter 5, sections 5.3 and 5.4]{G}, is at our disposition. 

The action of the distribution $\nabla K$
on the test function $\varphi$ is
\begin{equation*}
\begin{split}
\langle \nabla K, \varphi \rangle = & - \int K(x) \nabla \varphi(x)\, dx = -\lim_{\ep \rightarrow 0} 
\int_{|x|> \ep} K(x) \nabla \varphi(x)\, dx\\*[7pt]
& =  \lim_{\ep \rightarrow 0} \left(
\int_{|x|> \ep} \nabla K(x) \varphi(x)\, dx + \int_{|x|=\ep} K(x)\varphi(x) \frac{x}{|x|} d\sigma(x)\right)\\*[7pt]
& = \lim_{\ep \rightarrow 0} \left(
\int_{|x|> \ep} \nabla K(x) \varphi(x)\, dx + \varphi(0) \int_{|x|=\ep} K(x) \frac{x}{|x|} d\sigma(x)\right),
\end{split}
\end{equation*}
where $\sigma$ is surface measure.
If $K$ is even then $\int_{|x|=\ep} K(x) \frac{x}{|x|} d\sigma(x)=0$ and principal values of $\nabla K$ exist. Principal values
of $\nabla K$ also exist in the context of \cite{H} and \cite{ABC}, in which $K$ is homogeneous of order $-(N-1)$, because in this case
$\int_{|x|=\ep} K(x) \frac{x}{|x|} d\sigma(x)$ does not depend on $\ep.$ But for the kernel
$$
K(x) = \frac{\sin(\log(\frac{1}{|x|}))}{|x|^{N-1}} \frac{x_1}{|x|}, \quad 0 \neq x \in \RN,
$$
principal values of  $\nabla K$ do not exist, because $\lim_{\ep \rightarrow 0}\int_{|x|=\ep} K(x) \frac{x}{|x|} d\sigma(x)$ does not exist.

However, coming back to the general case, for some sequence $(\ep_j)_j$ tending to $0$ as $j$ tends to $\infty$ the limit
$$
\lim_{j \rightarrow \infty}\int_{|x|=\ep_j} K(x) \frac{x}{|x|} d\sigma(x) = \vec{L}
$$
does exist. This yields  existence of the principal value
$$
W(\varphi)  = \lim_{j \rightarrow \infty} \int_{|x|> \ep_j} \nabla K(x) \varphi(x)\, dx 
$$
and  the distributional identity
\begin{equation}\label{distr}
 \nabla K = W + \vec{L} \, \delta_0.
\end{equation}

Set $T(\varphi)= W\star \varphi,$ for $\varphi$ a test function. Then $T$ extends to a bounded operator on $L^p(\RN), \; 1< p< \infty$
 and to an operator that sends, with bounds, finite Radon measures into $L^{1,\infty}(\RN)$ \cite[Theorem 5.3.3] {G}. The result is indeed stated there only for $L^1(\RN)$ functions, 
 but it is well known that the Calder\'on-Zygmund decomposition and weak type estimates for singular integrals work also for general finite Radon measures. See, for instance, \cite[p.67]{H} for the Calder\'on-Zygmund decomposition and \cite[Theorem 20.26 and its proof]{M} for the weak type estimate for the singular integral. If one sets
$$
T_{\ep}(f)(x) = \int_{|y|> \ep} f(x-y) \nabla K(y)\,dy, \quad x \in \RN, 
$$
and
$$
T^*(f)(x)= \sup_{\ep >0}|T_{\ep}(f)(x)|, \quad x \in \RN,
$$
then the maximal singular integral $T^*$ is a bounded operator from $L^p(\RN), \; 1< p< \infty$ into itself and sends, with bounds, 
finite Radon measures into $L^{1,\infty}(\RN)$ \cite[Theorem 5.3.5 and Theorem 5.4.5]{G} for $L^1(\RN)$
functions and  \cite[Theorem 20.26]{M} for general finite Radon measures. See also \cite[Chapter 2]{T}, 
where the theory is developed 
even for a non-doubling underlying measure.

By \eqref{distr}, given a finite Radon measure $\mu,$ we have
\begin{equation}\label{gradK}
\nabla K \star \mu = T(\mu) + \vec{L} \,\mu.
\end{equation}
Recall that $T(\mu)$ and $\vec{L}$ depend on the sequence $(\ep_j)_j$ but the left hand side above
does not.
Now we want to assign a precise value to $\nabla K \star \mu$ at almost all points $a \in \RN.$ By the weak type
estimate for the maximal singular integral and the existence of principal values of $\nabla K$ on test functions 
along the sequence $(\ep_j)_j,$ one obtains, by standard reasoning,  that the principal value
\begin{equation}\label{pv}
 \lim_{j \rightarrow \infty} \int_{|x|> \ep_j} \nabla K(a-y) \, d\mu(y)
\end{equation}
exists for almost all $a. $ We then define  $T(\mu)(a)$ as this principal value. 
We define the value of $\mu$ at $a$  as  
$$\tilde{\mu}(a) = \lim_{r \rightarrow 0} \frac{\mu B(a,r)}{|B(a,r)|},$$
whenever this limit exists. We know that the limit above exists a.e. and coincides with  the absolutely continuous part of $\mu.$

Set
$$
(\nabla K \star \mu)(a)= T(\mu)(a) + \vec{L} \,\tilde{\mu}(a)
$$
at those points $a$ where the principal value \eqref{pv} and the limit in the definition of $\tilde{\mu}$  exist.
Note that this definition is just the identity 
$$
\nabla K  \star \varphi = T(\varphi)+ \vec{L} \, \varphi
$$
for $\mu = \varphi(x) dx, \; \varphi$ a test function.  Moreover, one has the weak type inequality
$$
t\, \left| \{a : \left|(\nabla K\star \mu)(a)\right| > t \}\right| \le C \, \|\mu  \|,
$$
because $\left|(\nabla K\star \mu)(a)\right| \le |T(\mu)(a)|+ |\vec{L}| \, M(\mu)(a),$ where $M(\mu)$ is the 
Hardy-Littlewood maximal  operator of $\mu.$
We claim that for $u=K\star\mu$ and for almost all $a \in \RN$
 $$
  \lim_{r \rightarrow 0} \frac{1}{\operatorname{Cap}(B(a,r))} \sup_{t>0} t \, \operatorname{Cap}\left(\{x\in B(a,r) : 
 \frac{|u(x)-u(a) -(\nabla K \star \mu)(a)\cdot(x-a) |}{|x-a|} > t \}\right) =0,
 $$
 which proves Theorem A with an explicit expression for the gradient of $u.$
This is clearly true for measures $\mu$ of the type
$$
\varphi(x) \,dx + \nu
$$
 with $\varphi \in C^\infty_0(\RN)$ and $\nu$ a singular measure supported on a closed set of Lebesgue measure zero. 
 Since this set of measures is dense in the total variation norm in the set of all finite Radon measures, we only need to prove
 that
 \begin{equation}
 \begin{split}
 & S(\mu)(a):=\\*[7pt] & \sup_{r>0}  \frac{1}{\operatorname{Cap}(B(a,r))} \sup_{t>0} t \, \operatorname{Cap} \left(\{x\in B(a,r) : 
 \frac{|u(x)-u(a) -(\nabla K \star \mu)(a)\cdot(x-a) |}{|x-a|} > t \}\right)
 \end{split}
 \end{equation}
 satisfies the weak type estimate 
 $$
 t\, \left| \{x : S(\mu)(x) > t \}\right| \le C \, \|\mu  \|.
 $$
To prove the preceding inequality assume without loss of generality  that $a=0$ and that $\mu$ is a positive measure. Inspired by \cite[Chapter 8, section 1]{S},  we proceed as follows :
\begin{equation*}
\begin{split}
\frac{1}{|x|} \left|u(x)-u(0) -(\nabla K \star \mu)(0)\cdot x\right| & \le 
\frac{1}{|x|} \left|u(x)-u(0) -(\nabla K \star \chi_{B(0,2 |x|)^c}\mu)(0)\cdot x \right| \\*[7pt]
& +  \left| (\nabla K \star \mu) (0) - (\nabla K \star \chi_{B(0,2|x|)^c}\mu)(0) \right| \\*[7pt]
= I +II.
\end{split}
\end{equation*}
The term $II$ is estimated by $|(\nabla K\star\mu)(0) | + |T^*(\mu)(0) |$, which is good because both terms have the weak type estimate
(since principal values do not exist in general one cannot say that $II$ tends to $0$ with $|x|$). The term $I$ is less than or equal to a sum of three terms $A+B+C$, where
$$
A= \frac{1}{|x|} \left| \int_{|y|> 2 |x|} K(x-y) \,d\mu(y) - \int_{|y|> 2 |x|} K(-y) \,d\mu(y)- \int_{|y|> 2 |x|}
\nabla K(-y) \,d\mu(y) \cdot x  \right|,
$$
$$
B=\frac{1}{|x|} \left| \int_{|y|\le 2 |x|} K(x-y) \,d\mu(y) \right|
$$
and
$$
C=\frac{1}{|x|} \left| \int_{|y|\le 2 |x|} K(-y) \,d\mu(y) \right|.
$$
The integrand in $A$ is estimated by the Taylor remainder of order $2$ and we get
$$
A \le C\, |x| \int_{|y|> 2 |x|} \frac{1}{|y|^{N+1}}\, d\mu(y) \le C \, M(\mu)(0).
$$
 This is a well known elementary inequality, which follows 
readily by splitting the domain of integration in dyadic annuli.
The term $C$ is no larger than 
$$
 \frac{1}{|x|}  \int_{|y|\le 2 |x|} \frac{1}{|y|^{N-1}}\, d\mu(y) \le C \, M(\mu)(0),
$$
as one can see easily, again splitting the domain of integration in dyadic annuli.

For the term $B=B(x)$ we prove the weak type capacitary estimate
$$
 \sup_{r>0}  \frac{1}{\operatorname{Cap}(B(0,r))} \sup_{t>0} t \, \operatorname{Cap}\left(\{x\in B(0,r) : 
 B(x) > t \}\right) \le C \, M(\mu)(0).
$$
Observe that if $|x| < r, $ then
$$
B(x) \le \frac{C}{|x|}  \int_{|y|\le 2 |x|} \frac{1}{|x-y|^{N-1}}\, d\mu(y) \le C\, \int_{|y|\le 2 r} \frac{1}{|x-y|^{N-1}}\,
\frac{d\mu(y)}{|y|}.
$$
By the weak capacitary estimate \eqref{wcap} (recall that ${\rm{Cap}}(B(a,r)) = c\, r^{N-1}$)
$$
\sup_{r>0}  \frac{1}{\operatorname{Cap}(B(a,r))} \sup_{t>0} t \, \operatorname{Cap}\left(\{x\in B(0,r) : 
 B(x) > t \}\right) \le \frac{C}{r^{N-1}} \, \int_{|y|\le 2 r} \frac{d\mu(y)}{|y|} \le C\, M(\mu)(0),
$$
where the last inequality follows readily by the usual method of decomposing the domain of integration in dyadic annuli ($N\ge 2$).
The proof is now complete.
 
The point of the proof is to avoid appealing to the Calder\'on-Zygmund decomposition of the measure $\mu$,
which has already been used in proving
the weak type inequality for $T^*(\mu).$ In particular this argument also proves readily the result in \cite{ABC} dealing with 
homogeneous kernels $K$,  for which principal values exist.

\section{Proof of Theorem B}
Given $x$ and $y$ in $\RN,$ set $r=|x-y|>0.$  Split the measure $\mu$ as $\mu =  \chi_{B(x,2r)} \mu + \chi_{B(x,2r)^c} \mu.$ Let us first take care of $\nu = \chi_{B(x,2r)^c} \mu.$ We have
\begin{equation*}
\begin{split}
& \left(K\star \nu \right)(y)-\left(K\star \nu \right)(x) = \int_{B(x,2r)^c}  \left(K(y-w) -K(x-w)\right) \,d\mu(w) \\*[7pt]
 &= \int_{B(x,2r)^c}  \left(K(y-w) -K(x-w)- \nabla K(x-w)\cdot (y-x)\right) \,d\mu(w) \\*[7pt] &  + \int_{B(x,2r)^c} 
 \left( \nabla K(x-w)\cdot (y-x)\right) \,d\mu(w) = I+II.
\end{split}
\end{equation*}
Estimate the integrand in $I$ by the Taylor remainder of order $2$  to get
\begin{equation*}
\begin{split}
| I | & \le  \int_{B(x,2r)^c}  \, \sup_{ |\xi| \ge |w-x| /2} \left| \nabla^2 K(\xi) \right| |x-y|^2 \,d\mu(w) \\*[7pt]
 & \le C\, |x-y| \,r \int_{B(x,2r)^c} \frac{d\mu(w)}{|w-x|^{N+1}} \\*[7pt] & \le C\, |x-y| \,M(\mu)(x).
\end{split}
\end{equation*}

For the term $II$  we obtain
$$
|II| = | \int_{B(x,2r)^c} \nabla K(x-w) \, d\mu(w) \cdot (y-x)| \le |x-y| \, T^*(\mu) (x).
$$

From now on we assume that $\mu$ lives in $B(x,2r).$   Let $\varphi_r(z) = \frac{1}{r^N} \varphi(\frac{z}{r} ), $  $\varphi$ 
an even continuously differentiable function on $\RN$ with compact support in the unit ball and integral equal to $1.$  
The required inequality clearly follows by combining the two estimates 
\begin{equation}\label{reg}
\left| \left(K\star \mu \right)(x)-\left(K\star \mu \star \varphi_r \right)(x) \right|  \le  C\, |x-y| \,M(\mu)(x)
\end{equation}
and
\begin{equation}\label{comp}
\left|\left(K\star \mu \star \varphi_r \right)(x) -\left(K\star \mu \star \varphi_r \right)(y) \right| \le C\,|x-y|\,M(\mu)(x).
\end{equation}
Note that by symmetry \eqref{reg} holds also with $x$ replaced by $y.$

We start by proving \eqref{comp}, which is straightforward :
\begin{equation*}
\begin{split}
\left|\left(K\star \mu \star \varphi_r \right)(x) -\left(K\star \mu \star \varphi_r \right)(y) \right| &= 
|\int \left(\varphi_r(x-z) -\varphi_r(y-z)\right) (K\star \mu)(z) \,dz | \\*[7pt] &  \le C\, 
\int_{|z-x|< 2r} \	\frac{|x-y|}{r^{N+1}} \, |(K\star \mu)(z)|\,dz \\*[7pt] & 
\le C \,\frac{|x-y|}{r^{N+1}} \,\int_{|w-x|<2r} \int_{|z-x|< 2r} \frac{dz}{|z-w|^{N-1}} \, d\mu(w) \\*[7pt] & 
\le C\, \frac{|x-y|}{r^{N}} \mu(B(x,2r)) \\*[7pt] &  \le C\, |x-y|\,M(\mu)(x).
\end{split}
\end{equation*}

We turn now to the proof of \eqref{reg}. We have, assuming that $\mu$ is absolutely continuous with a smooth density, 
(we will discuss later how to avoid this extra hypothesis),
\begin{equation*}
\begin{split}
D : &=   \left(K \star \mu \right)(x)-\left(K\star \mu \star \varphi_r \right)(x)  = 
\int \left( \left(K \star \mu \right)(x) - \left( K\star \mu \right)(z) \right) \varphi_r (x-z) \,dz   \\*[7pt]
 &= -\int  \left(\int_0^1   \left(\nabla K \star \mu \right)(x+t(z-x))  \cdot (z-x) \,dt  \right) \,   \varphi_r (x-z) \,dz  \\*[7pt] &  
 = -\int_0^1 \left(  \int   \left(\nabla K \star \mu \right)(\xi) \cdot \frac{\xi-x}{t}  \,  \varphi_r (\frac{x-\xi}{t}) 
 \,\frac{d\xi}{t^N}      \right) \,dt \\*[7pt] & 
 = -r\,\int_0^1 \left(  \int  \left(\nabla K \star \mu \right)(\xi) \cdot  \psi_R (x-\xi) \, d\xi      \right) \,dt,
\end{split}
\end{equation*}
where $R= t r$  and  $\psi(z) = z \varphi(z).$  Note that $\psi$ is an odd function and thus has zero integral 
(this will be used later on). 
By \eqref{gradK} 
\begin{equation*}
\begin{split}
-D &  = r\, \int_0^1 \left(  \int  T (\mu)(\xi) \cdot \psi_R(\xi-x) \, d\xi      \right) \,dt + r\, \int_0^1 \left(  \int  \vec{L} \cdot \psi_R(\xi-x) \, d\mu(\xi)      \right) \,dt \\*[7pt] & 
= I+II.
\end{split}
\end{equation*}
The second term is directly estimated by
$$
|II| \le C\, r \,\frac{\mu(B(x,R))}{R^N} \le C\, |x-y|\, M(\mu)(x).
$$

Let $S$ be the adjoint of $T.$  Then the term $I$ is 
\begin{equation}\label{atom}
\begin{split}
I &  =  r\, \int_0^1 \left(  \int  S \left(  \psi_R (\xi-x)\right) (w)\, d\mu(w)  \right) \,dt.
\end{split}
\end{equation}

 We claim that
\begin{equation}\label{claim}
\left| S\left(  \psi_R (\xi-x) \right) (w)  \right| \le C\,\left( \frac{1}{R^N}  \chi_{B(x,2R)} (w)+ R \, \frac{1}{|w-x|^{N+1}} \,\chi_{B(x,2R)^c} (w)  \right), \quad w \in \mathbb{R}^N.
\end{equation}
From this it follows readily that
$$
|I| \le C\,  |x-y|\, M(\mu)(x).
$$
To prove \eqref{claim} note that the operator $S$ acts as follows
$$
 S \left(  \psi_R (\xi-x) \right) (w) = 	\lim_{j \rightarrow \infty} \int_{|\xi-w|>\epsilon_j} \nabla K(\xi-w) \psi_R(\xi-x)\,d\xi.
$$
Assume without loss of generality that $x=0$ and distinguish $3$ cases.

Case 1 : $2R < |w|.$ Then, since $\int \psi_R(\xi)\,d\xi =0,$
$$
S \left(  \psi_R (\xi) \right) (w) = \int  \nabla K(\xi-w) \psi_R(\xi)\,d\xi = \int  \left(\nabla K(\xi-w)-\nabla K(-w)\right) \psi_R(\xi)\,d\xi 
$$
and so
$$
| S \left(  \psi_R (\xi) \right) (w)| \le C\, \int  \frac{|\xi|}{|w|^{N+1}} |\psi_R(\xi) |\,d\xi \le C\, \frac{R}{|w|^{N+1}}.
$$

Case 2 : $R < |w| < 2R.$  We have
\begin{equation*}
\begin{split}
S \left(  \psi_R (\xi) \right) (w) & = \int  \nabla K(\xi-w) \psi_R(\xi)\,d\xi = - \int  K(\xi-w) 	\nabla \psi_R(\xi)\,d\xi 
\end{split}
\end{equation*}
and hence
\begin{equation*}
\begin{split}
\left|S \left(  \psi_R (\xi) \right) (w) \right |& \le  \int  |K(\xi-w) |	|\nabla \psi_R(\xi)|\,d\xi \\*[7pt] & \le C\, \int_{|\xi|< R} \frac{1}{|\xi-w |^{N-1}} \,\frac{1}{R^{N+1}}\,d\xi \\*[7pt] & \le C\,  \frac{1}{R^{N}},
\end{split}
\end{equation*}
where we have applied a particular case of the general simple inequality
$$
\int_{E}  \frac{1}{|\xi-w |^{N-1}} \,d\xi \le C\, |E|^{\frac{1}{N}}, \quad E \;\; \text{measurable}\;\;\subset \RN.
$$

Case 3 : $|w| < R.$ This is as in case $2$, except that an additional boundary term appears, owing to the singularity at $w$ :
\begin{equation*}
\begin{split}
S \left(  \psi_R (\xi) \right) (w) &=  \lim_{j \rightarrow \infty} \int_{|\xi-w|> \epsilon_j}
\nabla K(\xi-w) 	
\psi_R(\xi)\,d\xi \\*[7pt] & = - \lim_{j \rightarrow \infty} \int_{|\xi-w|> \epsilon_j}  K(\xi-w) 
\nabla \psi_R(\xi)\,d\xi -  \int_{|\xi-w|=\epsilon_j}   K(\xi-w) \,\psi_R(\xi)\, \frac{\xi-w}{|\xi-w|}\,d\sigma(\xi) 
 \\*[7pt] & 
\end{split}
\end{equation*}
and thus
\begin{equation*}
\begin{split}
\left|S \left(  \psi_R (\xi) \right) (w) \right |& \le   C\, \int_{|\xi|< R} \frac{1}{|\xi-w |^{N-1}} \,\frac{1}{R^{N+1}}\,d\xi
+ C\,\limsup_{j\rightarrow \infty} \int_{|\xi-w|=\epsilon_j} \frac{1}{{\epsilon_j}^{N-1}} \,  \frac{1}{R^{N}}\,d\sigma(\xi)\\*[7pt] & \le C\,  \frac{1}{R^{N}},
\end{split}
\end{equation*}
which completes the proof of the claim \eqref{claim} and of \eqref{reg}.

Some words are in order to explain how we get rid of the smoothness hypothesis on $\mu.$
We used smoothness to prove \eqref{reg}. Regularize $\mu$ by convolving with $\varphi_{\frac{1}{n}}, \; n=1,2, \dots$ We get smooth functions 
$\mu_n = \mu \star \varphi_{\frac{1}{n}}.$ Then \eqref{reg} holds with $\mu$ replaced by $\mu_n.$ Now $K\star\mu_n$ converges 
to $K\star \mu$ in $L_{\text{loc}}^p(\RN), \; 1\le p <\frac{N}{N-1},$ hence a.e., passing to a subsequence if necessary. 
Therefore convergence in the left hand sides is a.e., as desired.  For the right hand sides just observe that 
$$
M(\mu_n)(x) \le C\, M(\mu)(x), \quad x \in \RN.
$$
The proof is now complete.

\subsection*{Remark} If the kernel $K$ is homogeneous of degree $-(N-1)$ the estimate of the term \eqref{atom} is much more direct.
On one hand,  $S(\psi_R) = S(\psi)_R,$  by homogeneity. On the other hand, since $\psi$ has zero integral, 
$S(\psi)$ is a $C^\infty(\RN)$ function such that  $ |S(\psi)(x) |	\le C\, |x|^{-(N+1)}, $  for large $|x|.$ 
Thus  the least radial majorant of $|S(\psi)|$ is a continuous function in $L^1(\RN)$ and so, by the general theory of the maximal operator, we have \cite[2.1.12, p.92]{G}
$$
\int \left| S \left(  \psi_R (x-\xi) \right) (w)\right|  \, d\mu(w)  \le C\, M(\mu)(x), \quad x \in \RN.
$$

\section{Level sets of newtonian potentials}
Let $P = |x|^{-(N-2)} \star \mu$ be the newtonian potential of a finite Radon measure $\mu,$  $N \ge 3.$ 
In the plane $P$ is the logarithmic potential $\log(|x|^{-1}) \star \mu$ and $\mu$ is a compactly supported Radon measure. 
The function $P$ is defined except for a set of zero newtonian capacity (logarithmic capacity in the plane). For a given constant $c$ set
$$
E= \{x \in \RN : P(x) =c \}
$$
In \cite{APR} one proves that the absolutely continuous part of $\mu$ vanishes on $E.$ See the introduction in \cite{APR} for the origin of this question. In \cite{CV} this was proved for the equilibrium measure $\mu$  of a compact set $E,$  for which $P$ is constant on $E$ except for a set of zero newtonian  capacity (logarithmic capacity in the plane). The result can be rephrased in the plane by saying that  harmonic measure is singular with respect to area, an old result of Oksendal \cite{O} from the seventies. 

The argument goes as follows. The gradient of $P$ is given by 
$$
\nabla P =  c_N \,\frac{x}{|x|^{N}} \star \mu,
$$
with $c_N= -(N-2), \; N>2$ and $c_N = -1,  \; N=2.$
Hence $\nabla P$ is locally in $L^p(\RN), \; 1\le p < N/(N-1).$  Thus $P$ is differentiable in the $L^p$ sense a.e.
 and, in particular, approximately differentiable a.e. \cite[p.230 and p.233]{EG}. Therefore $\nabla P =0$ a.e. on each level set of $P$ \cite[p.232]{EG}. 
 
We have
$$
\partial_j P =  c_N \, x_j / |x|^{N} \star \mu = K_j \star \mu, \quad 1 \le j \le N,
$$ 
with $K_j$ a smooth homogeneous kernel of degree $-(N-1).$ By Theorem A the function $\partial_j P$  is differentiable in the $L^p$ sense a.e., $1\le p < N/(N-1)$ and therefore, as before,  $\partial_j^2 P (x) =0$ a.e. on $E,$  because  $\partial_j P (x) =0, $ a.e. on $E.$ By the precise form of the gradient of $\partial_j P =K_j \star \mu$ given in the proof  of Theorem A 
 $$
 \partial_j^2 P (x) = c_N\, \left( \operatorname{p.v.} \frac{|x|^2-N x_j^2}{|x|^{N+2}} \star \mu (x)+V_N \, \tilde{\mu}(x) \right),\quad \text{for almost all} 	\;\;x,
$$
where $V_N$ is a dimensional constant (the volume of the unit ball in $\RN$) and $\tilde{\mu}(x)$ is the absolutely continuous part of $\mu$ at the point $x.$ Hence
$$
0 = \sum_{j=1}^N \partial_j^2 P (x) = c_N\, N V_N \,\tilde{\mu}(x), \quad \text{a.e. on }\;\; E,
$$
as desired.

A final remark is that using the results of \cite{CV} one can show that in the first step of the preceding argument one gets,
for $N\ge 3,$ that $\nabla u=0$ on $E$ except for a set of zero $C^1$ harmonic capacity, which is stronger by one dimension than saying
that $\nabla u=0$ a.e. on $E.$ This can be likely exploited to obtain sharp conclusions in situations similar to those 
envisaged in this section.

\begin{gracies}
The author would like to express his gratitude to Juli\`a Cuf\'i for many valuable suggestions and to Augusto Ponce 
for an illuminating correspondence that improved the exposition significantly.
This research was partially supported by the grants 2017SGR395
(Generalitat de Catalunya) and  MTM2016--75390 (Mi\-nis\-terio de
Educaci\'{o}n y Ciencia).
\end{gracies}

\begin{tabular}{l}
Joan Verdera\\
Departament de Matem\`{a}tiques\\
Universitat Aut\`{o}noma de Barcelona\\
08193 Bellaterra, Barcelona, Catalonia\\
{\it E-mail:} {\tt jvm@mat.uab.cat}
\end{tabular}


\end{document}